\documentclass[a4paper, 11pt]{article}
\usepackage{amsmath}
\usepackage{amsfonts}
\usepackage{amssymb}
\usepackage[english]{babel}
\usepackage{graphicx}
\usepackage{amsthm}
\usepackage[title]{appendix}

\usepackage{longtable}
\usepackage{tabularx}
\allowdisplaybreaks

\newtheorem{theorem}{Theorem}[section]

\newtheorem{lemma}[theorem]{Lemma}

\newtheorem{construction}[theorem]{Construction}
\def\whitebox{{\hbox{\hskip 1pt
 \vrule height 6pt depth 1.5pt
 \lower 1.5pt\vbox to 7.5pt{\hrule width
    3.2pt\vfill\hrule width 3.2pt}%
 \vrule height 6pt depth 1.5pt
 \hskip 1pt } }}
\def\qed{\ifhmode\allowbreak\else\nobreak\fi\hfill\quad\nobreak
     \whitebox\medbreak}

\newcommand{\ignore}[1]{}

\begin{document}
\baselineskip 16pt
\title{The Hamilton-Waterloo problem with $C_8$-factors and $C_m$-factors }

\author{\small  L. Wang and H. Cao \thanks{Research
supported by the National Natural Science Foundation of China
under Grant 11571179, the Natural Science Foundation of Jiangsu
Province under Grant No. BK20131393, and the Priority Academic
Program Development of Jiangsu Higher Education Institutions.
E-mail: {\sf caohaitao@njnu.edu.cn}} \\
\small Institute of Mathematics, \\ \small  Nanjing Normal
University, Nanjing 210023, China}

\date{}
\maketitle
\begin{abstract}

In this paper, we almost completely solve the Hamilton-Waterloo problem
 with $C_8$-factors and $C_m$-factors where the number of vertices is a multiple of $8m$.


\smallskip
\end{abstract}

\section{Introduction}

In this paper, we denote the cycle of length $k$ by
 $C_k$  and the complete $u$-partite graph with $u$ parts of
size $g$ by $K_u[g]$.
Let $H$ be a graph.
A {\it factor} of $H$ is a spanning subgraph of $H$.
Suppose $G$ is a subgraph of $H$, a {\it $G$-factor} of  $H$ is a set
 of edge-disjoint subgraphs of $H$,
each isomorphic to $G$. A {\it $G$-factorization} of $H$ is a set of edge-disjoint $G$-factors of $H$.
 A $C_k$-factorization of  $H$ is a partition of
$E(H)$ into $C_k$-factors.  Many papers introduced $C_k$-factorization of $K_u[g]$, see
\cite{ASSW,AH, HS, L1, L, LL,PWL, R}.

\begin{theorem}\label{Kv}
There exists a $C_k$-factorization of $K_u[g]$ if and only if $
g(u-1 )\equiv 0\pmod 2$, $gu\equiv 0\pmod k$, $k$ is even when
$u=2$, and $(k,u,g)\not\in\{(3,3,2),(3,6,2),(3,3,6),(6,2,6)\}$.
\end{theorem}

An {\it $r$-factor} is
a factor which is $r$-regular. Obviously, a 2-factor
consists of a collection of disjoint cycles.
 A {\it
$2$-factorization} of a graph $H$ is a partition of the edge set of $H$ into
2-factors.
The {\it Hamilton-Waterloo problem} is the problem of
determining whether $K_v$ (for $v$ odd) or $K_v$ minus a
$1$-factor (for $v$ even) has a $2$-factorization in which there
are exactly $\alpha$ $C_m$-factors and $\beta$ $C_n$-factors.
The authors {\rm \cite{WCC}} generalize this problem to an $r$-regular graph $H$, and
use HW$(H;m,n;\alpha,\beta)$ to denote a $2$-factorization of $H$ (for $r$ even) or $H$ minus a
$1$-factor (for $r$ odd)
in which there are exactly $\alpha$ $C_m$-factors and $\beta$
$C_n$-factors. Denote by
HWP$(H;m,n)$ the set of $(\alpha,\beta)$ for which an
HW$(H;m,n;\alpha,\beta)$ exists. When $H=K_v$, the two notations
HW$(H;m,n;\alpha,\beta)$ and  HW$(H;m,n)$ are written as
HW$(v;m,n;\alpha,\beta)$ and  HW$(v;m,n)$, respectively.
It's easy to see that the necessary conditions for the existence
of an HW$(v;m,n;\alpha,\beta)$ are $m | v$ when $ \alpha >0$, $n |
v$ when $ \beta >0$ and $\alpha+\beta=\lfloor \frac{v-1}{2}
\rfloor$. When $\alpha\beta=0$, the existence of an
HW$(v;m,n;\alpha,\beta)$ has been completely solved, see Theorem~\ref{Kv}.

Many authors have considered the Hamilton-Waterloo problem for small values of $m$ and $n$.
A complete solution for the existence of an
HW$(v;3,n;\alpha,\beta)$ in the cases $n\in\{4,5,7\}$ is given in
\cite{ABBE,DQS,LF,OO,WCC}.
For the case
$(m,n)\in\{(3,15),(5,15),(4,6),(4,8),(4,16),(8,16)\}$, see
\cite{ABBE}. The existence of an HW$(v;4,n;\alpha,\beta)$ for odd
$n\geq3 $ has been solved except possibly when $v=8n$ and
$\alpha=2$, see \cite{KO,OO,WCC}.
The authors \cite{FH} give a complete solution for the
existence of an HW$(v;4,n;\alpha,\beta)$ for even $n \geq 4$.
It is shown in \cite{K} that the
necessary conditions for the existence of an
HW$(v;3,9;\alpha,\beta)$ are also sufficient except possibly when
$\beta=1$. Many infinite classes of HW$(v;3,3x;\alpha,\beta)$s are
constructed in \cite{AKK}.
For more results on the Hamilton-Waterloo problem, the reader can see \cite{BDT,DL2,HNR,LFS,LS,WLC}.
In this paper, we focus on the existence of an
HW$(8mt;8,m;\alpha,\beta)$.

\begin{theorem}
\label{8m}  For $ m \geq 3$ and $ t \geq 1$, $(\alpha,\beta)\in$ \rm{HWP}$ (8mt;8,m)$
if and only if  $\alpha, \beta \geq 0$ and $\alpha+\beta=4mt-1$,
except possibly when $\alpha \in \{1,2\}$ and $mt$ is odd or  $\alpha \in \{1,2,4,5,6\}$ and $mt$ is even.
\end{theorem}

\section{Decompositions of Cayley graphs }

We start with some definitions in graph theory. For more general concepts of
graph theory, see \cite{W}.

Given a graph $G$, $G[n]$ is the
{\it lexicographic product} of $G$ with the empty graph on $n$ points.
Specifically, the point set is $\{x_i:x\in V(G), i\in Z_n\}$ and
$(x_i,y_j) \in E(G[n])$ if and only if $(x,y) \in E(G),i,j \in Z_n$. In
the following we will denote by $C_m[n]$ the lexicographic product
of $C_m$ with the empty graph on $n$ points.
We have the following known results.

\begin{theorem}{\rm (\cite{CNT,PWL})}
\label{Cm}  There exists a $C_{m}$-factorization of $C_m[n]$ for $m \geq 3$ and $n \geq 1$
except for $(m,n)=(3,6)$ and $(m,n) \in \{(l,2) \ | \ l \geq 3 $ is odd $\}$.
\end{theorem}

\begin{theorem}{\rm (\cite{LD})}
\label{Cmn}  There exists a $C_{mn}$-factorization of $C_m[n]$ for $m \geq 3$ and $n \geq 1$.
\end{theorem}

\begin{theorem}{\rm (\cite{OO})}
\label{Cm4}  There exists a $C_{4}$-factorization of $C_m[4]$ for $m \geq 3$.
\end{theorem}

Let $\Gamma$ be a finite additive group and let $S$ be a subset of
$\Gamma \backslash \{0\}$ such that the opposite of every element
of $S$ also belongs to $S$. The Cayley graph over $\Gamma$ with
connection set $S$, denoted by $Cay(\Gamma, S)$, is the graph with
vertex set $\Gamma$ and edge set $E(Cay(\Gamma, S))= \{(a,b)| a,b
\in \Gamma, a-b \in S\}$. It is quite obvious that $Cay(\Gamma,
S)=Cay(\Gamma, \pm S)$.

\begin{theorem}{\rm (\cite{WCC})}
\label{d=012} Let odd $m \geq 3$, $n \geq 3$.  If $a\in Z_n$, $| \pm \{ 0,a,2a\} |=5$ and $(i,m)=1$, then there is a $C_m$-factorization of
$Cay(Z_m \times Z_n,  \{ \pm i\}  \times  (\pm \{0, a,  2a\} ) )$.
\end{theorem}

\begin{lemma}
\label{d}  Let $ m \geq 3$, even $n\geq4$ and $0 < d < n$ be coprime to $n$. There exist
two $C_n$-factors which form a $C_n$-factorization of
$Cay(Z_m \times Z_n,  \{ \pm 1\}  \times  \{ \pm d\} )$.
\end{lemma}

\noindent {\it Proof:}
 Let $C_j=((a_{j0},b_{j0}),(a_{j1},b_{j1}),\cdots,(a_{j,n-1},b_{j,n-1}))$, $1\leq j\leq 2$, where

\hspace{2.1cm} $a_{10}=a_{21}=0; \ a_{20}=a_{11}=1; \ a_{jt}=a_{j,t-2}, \ 2 \leq t \leq n-1$, \\ \vspace{5pt}
\hspace{4.2cm} $b_{jt}=td\pmod n$, $0\le t\le n-1$.

Since $(d,n)=1$, we know that $b_{jt}$, $0\le t\le n-1$, are all
distinct modulo $n$. Then each $C_j$ will generate a $C_n$-factor
by $(+1\pmod{m},-)$. Thus we obtain the required two
$C_n$-factors which form  a $C_n$-factorization of $Cay(Z_m \times
Z_n, \{\pm 1\}\times\{ \pm d\})$.\qed

\begin{lemma}
\label{d=2}  Let $ m,n \geq 4$ be even integers and $1\leq d \leq n/2-1$. There exist
two $C_m$-factors which form a $C_m$-factorization of
$Cay(Z_m \times Z_n,  \{ \pm 1\}  \times  \{ \pm d\} )$.
\end{lemma}

\noindent {\it Proof:}
 Let $C_j=((0,b_{j0}),(1,b_{j1}),\cdots,(m-1,b_{j,m-1}))$, $1\leq j\leq 2$, where

\hspace{2.1cm} $b_{10}=b_{21}=0$; $b_{20}=b_{11}=d$; $b_{jt}=b_{j,t-2}$, $2 \leq t \leq m-1$.

Each $C_j$ can generate a $C_m$-factor
by $(-,+1\pmod{n})$. Then the conclusion follows. \qed

\begin{lemma}
\label{d=0,4}  Let even $n \geq 4$, $m\geq 3$ when $d=0$ or even $m\geq 4$ when $d=n/2$. There is
a $C_m$-factorization of $Cay(Z_m \times Z_n,  \{ \pm 1\}  \times  \{ d\} )$.
\end{lemma}

\noindent {\it Proof:} The cycle $(0_{0},1_{d},2_{0},3_{d},\dots ,(m-2)_{0},(m-1)_{d})$ can generate a $C_m$-factor
by $(-,+1\pmod{n})$ which is  a $C_m$-factorization of $Cay(Z_m \times
Z_n, \{\pm 1\}\times\{ d\})$.\qed

\begin{lemma}
\label{d=ab}  Let odd $ m \geq 3$, even $ n \geq 4$, and let $a,b \in Z_n$ with $| \pm \{ a,b,a+b\} |=6$. There exist
six $C_m$-factors which form a $C_m$-factorization of
$Cay(Z_m \times Z_n,  \{ \pm 1\}  \times (\pm \{ a,b,a+b\}) )$.
\end{lemma}

\noindent {\it Proof:}
 Let $C_j=((0,0),(1,b_{j1}),\cdots,(m-1,b_{j,m-1}))$, $1\leq j\leq 6$, where

\hspace{2.1cm} $b_{11}=-b_{41}=-b_{22}=b_{52}=a$; $b_{21}=-b_{51}=-b_{32}=b_{62}=b$;

\hspace{3.8cm} $b_{31}=-b_{61}=-b_{12}=b_{42}=-(a+b)$;  

\hspace{2.2cm} $b_{jt}=0$ when $t\geq 3$ is odd or $b_{jt}=b_{j1}$ when $t\geq 3$ is even.

Each $C_j$ can generate a $C_m$-factor
by $(-,+1\pmod{n})$. Then the conclusion follows.\qed

\begin{lemma}
\label{d=0d}  Let odd $ m \geq 3$, even $ n \geq 4$ and $1\leq d < n$. There exist three
 $C_m$-factors which form a $C_m$-factorization of
$Cay(Z_m \times Z_n,  \{ \pm 1\}  \times \{ 0,\pm d\} )$.
\end{lemma}

\noindent {\it Proof:}
 Let $C_j=((0,0),(1,b_{j1}),\cdots,(m-1,b_{j,m-1}))$, $1\leq j\leq 3$, where

\hspace{3.1cm} $b_{21}=-b_{31}=b_{12}=-b_{32}=d$; $b_{11}=b_{22}=0$;  
 
\hspace{2.1cm} $b_{jt}=0$ when $t\geq 3$ is odd or $b_{jt}=b_{j1}$ when $t\geq 3$ is even.

Each $C_j$ can generate a $C_m$-factor
by $(-,+1\pmod{n})$. We obtain the conclusion.\qed

\begin{lemma}
\label{d=41}  Let odd $ m \geq 5$. The graph $Cay(Z_m \times Z_8,  \{ \pm 1\}  \times  \{ \pm 1, 4\} )$
can be decomposed into two $C_8$-factors and a $C_m$-factor.
\end{lemma}

\noindent {\it Proof:}
For $m=5$,
the $C_5$-factor can be generated from $( 0_0, 1_1, 2_2, 3_3, 4_4)$ by $( -, +1\pmod 8)$.
Two $C_8$-factors are as below.

\vspace{5pt}

{\footnotesize
$\{( 0_0, 4_1, 0_2, 4_3, 0_4, 4_5, 0_6, 4_7), ( 0_{1+2i}, 1_{5+2i}, 2_{1+2i}, 3_{5+2i}, 4_{4+2i}, 3_{2i}, 2_{4+2i}$, $1_{2i})\ | \ 0\leq i \leq 3\}$,

$\{( 4_0, 0_1, 4_2, 0_3, 4_4, 0_5, 4_6, 0_7), ( 0_{2i}, 1_{4+2i}, 2_{3+2i}, 3_{2+2i}, 4_{1+2i}, 3_{5+2i}, 2_{6+2i}$, $1_{7+2i})\ | \ 0\leq i \leq 3 \}$.}

\vspace{5pt}

For $m\geq7$,
the $C_m$-factor can be obtained from $( 0_0, 1_1, 2_2, 3_3, 4_4, 5_0, 6_4,\cdots, i_0, (i+1)_4, \cdots,$ $ (m-2)_0, (m-1)_4)$ by $( -, +1\pmod 8)$.
Two $C_8$-factors are given as follows.

\vspace{5pt}

{\footnotesize
$\{( 0_0, (m-1)_1, 0_2, (m-1)_3, 0_4, (m-1)_5, 0_6, (m-1)_7), ( 0_{1+2j}, 1_{5+2j}, 2_{1+2j}$, $3_{5+2j}, 4_{4+2j}, 3_{2j}$, $2_{4+2j}, 1_{2j})$,
$( (i+1)_0, i_1, (i+1)_2, i_3, (i+1)_4, i_5, (i+1)_6, i_7) \ | \ 4\leq i \leq m-2, 0\leq j \leq 3 \}$,

$\{( (m-1)_0, 0_1, (m-1)_2, 0_3, (m-1)_4, 0_5, (m-1)_6, 0_7), ( 0_{2j}, 1_{4+2j}, 2_{3+2j}, 3_{2+2j}, 4_{1+2j}$, $3_{5+2j}$, $2_{6+2j}, 1_{7+2j})$,
$( i_0, (i+1)_1, i_2, (i+1)_3, i_4, (i+1)_5, i_6, (i+1)_7) \ | \ 4\leq i \leq m-2, 0\leq j \leq 3 \}$.}
\qed

For the following lemmas, we need a special 1-factorization of $K_8$ whose 7 1-factors are listed as below.

 \vspace{5pt}
{\footnotesize
$I_1=\{(0,1),(2,3),(4,5),(6,7)\}$,
$I_2=(\{(1,2),(3,4),(5,6),(7,0)\}$,
$I_3=\{(0,4),(1,5),(2,6),(3,7)\}$,

$I_4=\{(0,5),(1,7),(2,4),(3,6)\}$,
$I_5=\{(0,2),(1,3),(5,7),(4,6)\}$,
~$I_6=\{(0,6),(1,4),(2,7),(3,5)\}$,

$I_7=\{(0,3),(1,6),(2,5),(4,7)\}$.
}
Note that $I_i \cup I_{i+1}$ can form a 8-cycle for any $i \in \{1,3,5\}$.

\begin{lemma}
\label{d=2-13}  Let $ m \geq 3$. There exist
two $C_8$-factors which form a $C_8$-factorization of
$Cay(Z_m \times Z_8,  \{ \pm 1\}  \times  \{ 2\} )  \cup  mI_1 \cup mI_3$.
\end{lemma}

\noindent {\it Proof:}
Let $C_1=(0_0, 1_2, 1_6, 0_4, 0_5, 1_7, 1_3, 0_1)$,
$C_2=(1_0, 0_6, 0_7, 1_1, 1_5, 0_3, 0_2, 1_4)$.
Each $C_i$ will generate a $C_8$-factor by $(+1\pmod m, -)$.
Thus we get two $C_8$-factors which form a $C_8$-factorization of
$Cay(Z_m \times Z_8,  \{ \pm 1\}  \times  \{ 2\} )  \cup mI_1 \cup  mI_3$.\qed

\begin{lemma}
\label{d=6-27}  Let $ m \geq 3$. There exist
two $C_8$-factors which form a $C_8$-factorization of
$Cay(Z_m \times Z_8,  \{ \pm 1\}  \times  \{ 6\} )  \cup  mI_2 \cup  mI_7$.
\end{lemma}

\noindent {\it Proof:}
Let $C_1=(0_0, 1_6, 1_1, 0_3, 0_4, 1_2, 1_5, 0_7)$,
$C_2=(1_0, 0_2, 0_1, 1_7, 1_4, 0_6, 0_5, 1_3)$.
Each $C_i$ will generate a $C_8$-factor by $(+1\pmod m, -)$.
We get the conclusion.
\qed

\begin{lemma}
\label{d=4-456}  Let $ m \geq 3$.
The graph $Cay(Z_m \times Z_8,  \{ \pm 1\}  \times  \{ 4\} ) \cup m(\cup_{i=4}^6 I_i)$
can be decomposed into two $C_8$-factors and a 1-factor.
\end{lemma}

\noindent {\it Proof:}
The 1-factor is $\{(j_1, (j+1)_5), (j_2, (j+1)_6), (j_4, (j+1)_0), (j_7, (j+1)_3) \ | \ j \in Z_m\}$.
Let $C_1=(0_0, 0_2, 0_7, 0_5, 0_3, 0_1, 0_4, 0_6)$,
$C_2=(0_0, 0_5, 1_1, 1_7, 0_3, 0_6, 1_2, 1_4)$.
Each $C_i$ will generate a $C_8$-factor by $(+1\pmod m, -)$.
Thus, we obtain the required $C_8$-factors and 1-factor.
\qed

\begin{lemma}
\label{d=4+K8}  Let $ m \geq 3$.
The graph $Cay(Z_m \times Z_8,  \{ \pm 1\}  \times  \{ 4\} ) \cup   mK_8$
can be decomposed into four $C_8$-factors and a 1-factor.
\end{lemma}

\noindent {\it Proof:}
The 1-factor is $\{(j_3, (j+1)_7), (j_4, (j+1)_0), (j_5, (j+1)_1), (j_6, (j+1)_2)\} \ | \ j \in Z_m\}$. Let

 \vspace{5pt}
{\footnotesize
$C_1=(0_0, 0_2, 0_1, 0_3, 0_5, 0_4, 0_6, 0_7)$,
$C_2=(0_0, 0_3, 0_2, 0_4, 0_7, 0_5, 0_1, 0_6)$,

$C_3=(0_0, 0_4, 0_1, 0_7, 0_3, 0_6, 0_2, 0_5)$,
$C_4=(0_0, 0_1, 1_5, 1_6, 0_2, 0_7, 1_3, 1_4)$.}

For $1\leq i \leq 4$, each $C_i$ will generate a $C_8$-factor by $(+1\pmod m, -)$.
Thus we get the required conclusion.
\qed

 \section{Proof of Theorem~\ref{8m}}

In order to prove our main results, we still need the following construction.

\begin{construction} {\rm (\cite{WCC})}
\label{C-RGDD} If there exist an {\rm
HW}$( K_u[g] ; m ,n; \alpha,  \beta)$ and an {\rm
HW}$(g;m,n; \alpha', \beta' )$, then an {\rm
HW}$(gu;m,n;  \alpha + \alpha',  \beta +\beta' )$ exists.
\end{construction}

\begin{lemma}
\label{0248}  For $m\geq 3$ and $r \in \{0,2,4,8\}$, $(r,8-r)\in$ {\rm HWP}$(C_m[8];8,m)$.
\end{lemma}

\noindent {\it Proof:} We consider each of these cases in turn.

{\bf Case 1:} $r=0$.

The conclusion comes from Theorem~\ref{Cm}.

 \vspace{5pt}

{\bf Case 2:} $r=2$.

Two $C_8$-factors are given from a $C_8$-factorization of
$Cay(Z_m \times Z_8,  \{ \pm 1\}  \times  \{ \pm 3\} )$ by Lemma~\ref{d}.
The required six $C_m$-factors can be obtained as follows.

$(1)$ $m$ is even.

The graph $Cay(Z_m \times Z_8,  \{ \pm 1\}  \times  \{0,4 \} )$
can be decomposed into two $C_m$-factors by Lemma~\ref{d=0,4}. Similarly, $Cay(Z_m \times Z_8,  \{ \pm 1\}  \times  (\pm \{ 1,2\}) )$
can be partitioned into four $C_m$-factors from Lemma~\ref{d=2}.

$(2)$ $m$ is odd.

 Let $C_i^j=((0,b_{i0}^j),(1,b_{i1}^j), \cdots, (m-1,b_{i,m-1}^j))$, $1\leq i \leq6$, $j=1,2$, where

{\footnotesize \vspace{5pt}\noindent
\begin{tabular}{lllllll}
 \hspace{2.6cm} & $b_{10}^1=0$, & $b_{11}^1=0$, & $b_{12}^1=1;$  & $ b_{10}^2=1$, & $b_{11}^2=7$, & $b_{12}^2=0; $ \\
 \hspace{2.6cm} & $b_{20}^1=0$, & $b_{21}^1=2$, & $b_{22}^1=0;$  & $ b_{20}^2=1$, & $b_{21}^2=1$, & $b_{22}^2=3; $ \\
 \hspace{2.6cm} &  $b_{30}^1=0$, & $b_{31}^1=1$, & $b_{32}^1=7;$  & $ b_{30}^2=1$, & $b_{31}^2=0$, & $b_{32}^2=2; $ \\
 \hspace{2.6cm} &  $b_{40}^1=0$, & $b_{41}^1=6$, & $b_{42}^1=2;$  & $ b_{40}^2=1$, & $b_{41}^2=3$, & $b_{42}^2=7; $ \\
 \hspace{2.6cm} &  $b_{50}^1=0$, & $b_{51}^1=4$, & $b_{52}^1=4;$  & $ b_{50}^2=1$, & $b_{51}^2=5$, & $b_{52}^2=5; $ \\
 \hspace{2.6cm} &  $b_{60}^1=0$, & $b_{61}^1=7$, & $b_{62}^1=6;$  & $ b_{60}^2=7$, & $b_{61}^2=0$, & $b_{62}^2=7. $
\end{tabular}}

\vspace{5pt}

When $m\geq 5$, for any $t\geq 3$,
$(b_{1t}^1,b_{1t}^2,b_{2t}^1,b_{2t}^2,b_{3t}^1,b_{3t}^2,b_{4t}^1,b_{4t}^2,b_{5t}^1,b_{5t}^2,b_{6t}^1,b_{6t}^2)$=$(1,0,4,7,$ $1,4,0,5,5,6,5,6)$ for odd $t$
or $(1,0,0,3,7,2,2,7,4,5,6,7)$ for even $t$.
Each $F_i=\{C_i^j \ | \ j=1,2\}$ will generate a $C_m$-factor by $( -,+2\pmod 8)$.

 \vspace{5pt}

{\bf Case 3:} $r=4$.

Four $C_8$-factors come from a $C_8$-factorization of
$Cay(Z_m \times Z_8,  \{ \pm 1\}  \times (\pm \{ 1,3\}) )$ by Lemma~\ref{d}.
The graph $Cay(Z_m \times Z_8,  \{ \pm 1\}  \times  \{ 0\} )$
can be decomposed into a $C_m$-factor by Lemma~\ref{d=0,4}.
The other three $C_m$-factors are given as follows.

$(1)$ $m$ is even.

From Lemma~\ref{d=0,4}, $Cay(Z_m \times Z_8,  \{ \pm 1\}  \times  \{ 4\} )$
can be decomposed into a $C_m$-factor.
The graph $Cay(Z_m \times Z_8,  \{ \pm 1\}  \times  \{ \pm 2\} )$ will be partitioned into two $C_m$-factors from Lemma~\ref{d=2}.

$(2)$ $m$ is odd.

Let $C_i^j=((0,b_{i0}^j),(1,b_{i1}^j), \cdots, (m-1,b_{i,m-1}^j))$, $1\leq i \leq 3$, $1\leq j \leq 4$, where

{\footnotesize \vspace{5pt}\noindent
\begin{tabular}{lllllll}
 \hspace{2.6cm} & $b_{10}^1=0$, & $b_{11}^1=2$, & $b_{12}^1=4;$  & $ b_{10}^2=2$, & $b_{11}^2=0$, & $b_{12}^2=6; $ \\
 \hspace{2.6cm} & $b_{10}^3=1$, & $b_{11}^3=3$, & $b_{12}^3=5;$  & $ b_{10}^4=3$, & $b_{11}^4=1$, & $b_{12}^4=7; $ \\

 \hspace{2.6cm} & $b_{20}^1=0$, & $b_{21}^1=6$, & $b_{22}^1=2;$  & $ b_{20}^2=6$, & $b_{21}^2=0$, & $b_{22}^2=4; $ \\
 \hspace{2.6cm} & $b_{20}^3=1$, & $b_{21}^3=7$, & $b_{22}^3=3;$  & $ b_{20}^4=7$, & $b_{21}^4=1$, & $b_{22}^4=5; $ \\

 \hspace{2.6cm} &  $b_{30}^1=0$, & $b_{31}^1=4$, & $b_{32}^1=6;$  & $ b_{30}^2=6$, & $b_{31}^2=2$, & $b_{32}^2=0; $ \\
 \hspace{2.6cm} &  $b_{30}^3=1$, & $b_{31}^3=5$, & $b_{32}^3=7;$  & $ b_{30}^4=7$, & $b_{31}^4=3$, & $b_{32}^4=1. $
\end{tabular}}

\vspace{5pt}

For $m\geq 5$ and $t\geq 3$,
$(b_{1t}^1,b_{1t}^2,b_{1t}^3,b_{1t}^4,b_{2t}^1,b_{2t}^2,b_{2t}^3,b_{2t}^4,b_{3t}^1,b_{3t}^2,b_{3t}^3,b_{3t}^4)$=$(0,2,$ $1,3,4,6,$ $5,7,4,$ $6,5,7)$ for odd $t$
or $(4,6,5,7,2,4,3,5,6,0,7,1)$ for even $t$.
Each $F_i=\{C_i^j \ | \ 1\leq j \leq 4\}$ will generate a $C_m$-factor by $( -,+4\pmod 8)$.

 \vspace{5pt}

{\bf Case 4:} $r=8$.

The graph $C_m[4]$ can be decomposed into four $C_4$-factors for $m\geq 3$ by Theorem~\ref{Cm4}.
Give each vertex weight 2 to get four $C_4[2]$-factors.
From Theorem~\ref{Cmn}, each $C_4[2]$ of $C_4[2]$-factor can be partitioned into two $C_8$-factors.
Thus each $C_4[2]$-factors can be decomposed into two $C_8$-factors.
Finally, we get the required $C_8$-factors.
\qed

\begin{lemma}
\label{3-11} For $m\geq 3$ and $3\leq r\leq 11$, the graph $C_m[8] \cup  mK_8$
can be partitioned into $r$ $C_8$-factors, $11-r$ $C_m$-factors
and a 1-factor.
\end{lemma}

\noindent {\it Proof:} We distinguish six cases as below.

\vspace{5pt}

{\bf Case 1}: $r \in \{3,5,7,11\}$.

$C_m[8]$ can be partitioned into $r_1$ $C_8$-factors and $8-r_1$ $C_m$-factors
 for $m\geq 3$ and $r_1 \in \{0,2,4,8\}$ by Lemma~\ref{0248}.
Since $K_8$ can be decomposed into three $C_8$-factors
and a 1-factor by Theorem~\ref{Kv},
the graph $mK_8$ can be decomposed into three $C_8$-factors
and a 1-factor.
Thus we get $r=r_1+3$ $C_8$-factors, $8-r_1$ $C_m$-factors
and a 1-factor.

\vspace{5pt}

{\bf Case 2}: $r=4$.

$Cay(Z_m \times Z_8,  \{ \pm 1\}  \times  \{ 4\} ) \cup  mK_8$
can be decomposed into four $C_8$-factors and a 1-factor by Lemma~\ref{d=4+K8}.
A $C_m$-factor comes from a $C_m$-factorization of $Cay(Z_m \times Z_8,  \{ \pm 1\}  \times  \{ 0\} )$
by Lemma~\ref{d=0,4}.
$Cay(Z_m \times Z_8,  \{ \pm 1\}  \times  ( \pm \{ 1,2,3\}) )$
can be partitioned into six $C_m$-factors from Lemma~\ref{d=2} or Lemma~\ref{d=ab} when $m$ is even or odd.

\vspace{5pt}

{\bf Case 3}: $r=6$.

$Cay(Z_m \times Z_8,  \{ \pm 1\}  \times  \{ 4\} ) \cup mK_8$
can be decomposed into four $C_8$-factors and a 1-factor by Lemma~\ref{d=4+K8}.
The other two $C_8$-factors are given from a $C_8$-factorization of
$Cay(Z_m \times Z_8,  \{ \pm 1\}  \times \{ \pm 3\} )$ by Lemma~\ref{d}.
Five $C_m$-factors are as follows.

$(1)$ $m$ is even.

The graph $Cay(Z_m \times Z_8,  \{ \pm 1\}  \times  \{ 0\} )$
can be decomposed into a $C_m$-factor by Lemma~\ref{d=0,4}.
Further, $Cay(Z_m \times Z_8,  \{ \pm 1\}  \times ( \pm \{ 1,2\} ) )$
will be partitioned into four $C_m$-factors from Lemma~\ref{d=2}.

$(2)$ $m$ is odd.

By Lemma~\ref{d=012}, $Cay(Z_m \times Z_8,  \{ \pm 1\}  \times  (\pm \{ 0, 1, 2\} ))$
can be decomposed into five $C_m$-factors.

\vspace{5pt}

{\bf Case 4}: $r=8$.

From Lemma~\ref{d=4+K8}, $Cay(Z_m \times Z_8,  \{ \pm 1\}  \times  \{ 4\} ) \cup mK_8$
can be decomposed into four $C_8$-factors and a 1-factor.
The other four $C_8$-factors are given from a $C_8$-factorization of
$Cay(Z_m \times Z_8,  \{ \pm 1\}  \times  (\pm \{ 1,3\}) )$ by Lemma~\ref{d}.
Three $C_m$-factors are given as below.

$(1)$ $m$ is even.

The graphs $Cay(Z_m \times Z_8,  \{ \pm 1\}  \times  \{ 0\} )$ and
$Cay(Z_m \times Z_8,  \{ \pm 1\}  \times  \{ \pm 2\}  )$
can be decomposed into a $C_m$-factor and two $C_m$-factors by Lemmas~\ref{d=0,4} and \ref{d=2},
respectively.

$(2)$ $m$ is odd.

From Lemma~\ref{d=0d} $Cay(Z_m \times Z_8,  \{ \pm 1\}  \times  \{ 0, \pm 2\} )$
can be decomposed into three $C_m$-factors.

\vspace{5pt}

{\bf Case 5}: $r=9$.

The 1-factor is $mI_4$.

There exist
two $C_8$-factors from a $C_8$-factorization of
$Cay(Z_m \times Z_8,  \{ \pm 1\}  \times  \{ 2\} )  \cup mI_1 \cup mI_3$
by Lemma~\ref{d=2-13}.
The graph $Cay(Z_m \times Z_8,  \{ \pm 1\}  \times  \{ 6\} )  \cup mI_2 \cup mI_7$
can be decomposed into two $C_8$-factors by Lemma~\ref{d=6-27}.
The graph $\ m(I_5 \cup I_6)$
 can be decomposed into a $C_8$-factor since $I_5 \cup I_6$ can form a 8-cycle.
$Cay(Z_m \times Z_8,  \{ \pm 1\}  \times  \{ \pm 3\} )$
can be decomposed into two $C_8$-factors by Lemma~\ref{d}.
The other two $C_8$-factors and two $C_m$-factors are given as follows.

 $(1)$ $m$ is even.

By Lemma~\ref{d}, $Cay(Z_m \times Z_8,  \{ \pm 1\} \times \{ \pm 1\} )$
can be decomposed into two $C_8$-factors.
The graph $Cay(Z_m \times Z_8,  \{ \pm 1\}  \times  \{0,4\} )$
can be partitioned into two $C_m$-factors from Lemma~\ref{d=0,4}.

 $(2)$ $m$ is odd.

$\textcircled{\small{1}}$ $m=3$.

Two $C_3$-factors and two $C_8$-factors are as below.

 \vspace{5pt}
{\footnotesize
$\{( 0_0, 1_0, 2_1)$,
$( 0_1, 1_1, 2_0)$,
$( 0_2, 1_2, 2_3)$,
$( 0_3, 1_3, 2_2)$,
$( 0_4, 1_4, 2_5)$,
$( 0_5, 1_5, 2_4)$,
$( 0_6, 1_6, 2_7)$,
$( 0_7, 1_7, 2_6)\}$,

$\{( 0_0, 1_4, 2_0)$,
$( 0_4, 1_0, 2_4)$,
$( 0_1, 1_2, 2_1)$,
$( 0_2, 1_1, 2_2)$,
$( 0_3, 1_7, 2_3)$,
$( 0_7, 1_3$, $2_7)$,
$( 0_5, 1_6, 2_5)$,
$( 0_6, 1_5, 2_6)\}$.

$\{( 0_0, 1_1, 2_1, 0_5, 2_6, 1_2, 0_6, 1_7)$,
$( 1_0, 2_0, 0_7, 2_3, 1_4, 2_4, 0_3, 2_7)$,
$( 0_1, 2_2, 1_6$, $0_2$, $1_3, 0_4, 1_5, 2_5)\}$,

$\{( 0_0, 2_4, 1_3, 2_3, 0_4, 2_0, 1_7, 2_7)$,
$( 1_0, 0_1, 1_5, 2_1, 0_2, 2_6, 1_6, 0_7)$,
$( 1_1, 0_5, 1_4, 0_3, 1_2$, $2_2$, $0_6, 2_5)\}$.}

$\textcircled{\small{2}}$ $m \geq 5$.
%\Large{\textcircled{\small{2}}}

The graph $Cay(Z_m \times Z_8,  \{ \pm 1\}  \times  \{ \pm 1, 4\} )$
can be decomposed into two $C_8$-factors and a $C_m$-factor by Lemma~\ref{d=41}.
The last $C_m$-factor comes from a $C_m$-factorization of
$Cay(Z_m \times Z_8,  \{ \pm 1\}  \times  \{0\} )$ from Lemma~\ref{d=0,4}.

\vspace{5pt}

{\bf Case 6}: $r=10$.

Four $C_8$-factors are given from a $C_8$-factorization of
$Cay(Z_m \times Z_8,  \{ \pm 1\}  \times (\pm \{ 1,3\}) )$ by Lemma~\ref{d}.
The graph $Cay(Z_m \times Z_8,  \{ \pm 1\}  \times  \{ 2,4,6\} ) \cup mK_8$
can be decomposed into six $C_8$-factors and a 1-factor by Lemmas~\ref{d=2-13}-\ref{d=4-456}.
Further, the graph $Cay(Z_m \times Z_8,  \{ \pm 1\}  \times  \{ 0\} )$
can be partitioned into a $C_m$-factor by Lemma~\ref{d=0,4}.
\qed

Finally, we prove our main theorem.

 \vspace{5pt}

\noindent {\bf  Proof of Theorem~\ref{8m}:}  We distinguish two cases to discuss.

{\bf Case 1}: $mt$ is odd.

We start with an
HW$(K_{mt}[1];8,m;0,\frac{mt-1}{2})$ on the vertex set $Z_{mt} $ from Theorem~\ref{Kv}.
Give each vertex weight 8 to get $\frac{mt-1}{2}$ $C_m[8]$-factors which are denoted by $P_i$, $1\leq i \leq \frac{mt-1}{2}$ and $mtK_8$.
Each $P_i$ have $t$ $C_m[8]$s, denoted by $Q_{ij}$, $1\leq j \leq t$.

For $P_i$, $1\leq i \leq \frac{mt-3}{2}$, we only choose $x$ of them with $0\leq x \leq \frac{mt-3}{2}$ and replace each $C_{m}[8]$ in $x$ $C_m[8]$-factors with an HW$(C_{m}[8];8,m;8,0)$
from Lemma~\ref{0248}. Further, replace each $C_{m}[8]$ in the rest of $C_m[8]$-factors an HW$(C_{m}[8];8,m;0,8)$ from Lemma~\ref{0248}.

For $P_{\frac{mt-1}{2}}$ and $1\leq j \leq t$, the graph $Q_{\frac{mt-1}{2},j} \cup mK_8$
can be partitioned into $r$ $C_8$-factors, $11-r$ $C_m$-factors
and a 1-factor for $3\leq r\leq 11$ by Lemma~\ref{3-11}.
Put them together to get $r$ $C_8$-factors, $11-r$ $C_m$-factors
and a 1-factor on the vertex set $Z_{mt} \times Z_8$.

Hence, we obtain $\alpha=8x+r$ $C_8$-factors, $\beta=8\cdot (\frac{mt-3}{2}-x)+11-r=4mt-1-(8x+r)$ $C_m$-factors
and a 1-factor for $0\leq x \leq \frac{mt-3}{2}$ and $3\leq r\leq 11$. It's easy to check that $\alpha$ cover the integers from 3 to $4mt-1$.

{\bf Case 2}: $mt$ is even.

$(1)$ $\alpha =3$.

We can get the conclusion  by using Construction~\ref{C-RGDD} with an HW$(8;$ $
8, m; 3, 0)$ and an HW$( K_{mt}[8]; 8, m; 0,  4mt-4)$ from
Theorem~\ref{Kv}.

$(2)$ $\alpha \geq 7$.

Beginning with an
HW$(K_{mt}[1];8,m;0,\frac{mt-2}{2})$ with $\frac{mt-2}{2}$ $C_m$-factors and a 1-factor from Theorem~\ref{Kv} and
giving each vertex weight 8, we get $\frac{mt-2}{2}$ $C_m[8]$-factors which are denoted by $P_i$, $1\leq i \leq \frac{mt-2}{2}$, $mtK_8$ and $\frac{mt}{2}K_{2}[8]$.

For $P_i$, $1\leq i \leq \frac{mt-4}{2}$, we also choose $x$ of them with $0\leq x \leq \frac{mt-4}{2}$ and replace each $C_{m}[8]$ in $x$ $C_m[8]$-factors with an HW$(C_{m}[8];8,m;8,0)$
from Lemma~\ref{0248}. Next, replace each $C_{m}[8]$ in the rest of $C_m[8]$-factors an HW$(C_{m}[8];8,m;0,8)$ from Lemma~\ref{0248}.

It is similar to the first case, $P_{\frac{mt-2}{2}} \cup mtK_8$ can be decomposed into
 $r$ $C_8$-factors, $11-r$ $C_m$-factors
and a 1-factor on the vertex set $Z_{mt} \times Z_8$ for $m\geq 3$ and $3\leq r\leq 11$.

$\frac{mt}{2}K_{2}[8]$
can be partitioned into 4 $C_8$-factors by Lemma~\ref{Kv}.

We together get $\alpha=8x+r+4$ $C_8$-factors, $\beta=8\cdot (\frac{mt-4}{2}-x)+11-r=4mt-1-(8x+r+4)$ $C_m$-factors
and a 1-factor for $0\leq x \leq \frac{mt-4}{2}$ and $3\leq r\leq 11$.  We can check that $\alpha$ cover the integers from 7 to $4mt-1$.
\qed

\end{document}